\documentclass[a4paper,english,twoside,12pt]{smfart}

\usepackage{smfenum}
\usepackage{smfthm}
\usepackage{amssymb}
\usepackage{mathrsfs}

\def\RR{\mathbf{R}}
\def\ZZ{\mathbf{Z}} 

\def\A{{\rm A}}
\def\B{{\rm B}}
\def\C{{\rm C}}
\def\D{{\rm D}}
\def\E{{\rm E}}
\def\F{{\rm F}}
\def\G{{\rm G}}
\def\H{{\rm H}}
\def\I{{\rm I}}
\def\J{{\rm J}}
\def\K{{\rm K}}
\def\L{{\rm L}}
\def\M{{\rm M}}
\def\N{{\rm N}}
\def\P{{\rm P}}
\def\Q{{\rm Q}}

\def\SS{{\rm S}}
\def\T{{\rm T}}
\def\U{{\rm U}}
\def\V{{\rm V}}
\def\W{{\rm W}}
\def\X{{\rm X}}
\def\Y{{\rm Y}}
\def\Z{{\rm Z}}

\def\Aa{\mathscr{A}}
\def\Bb{\mathscr{B}}

\def\Dd{\mathscr{D}}
\def\Ee{\mathscr{E}}
\def\Ff{\mathscr{F}}

\def\Hh{\mathscr{H}}

\def\Mm{\mathscr{M}}
\def\Nn{\mathscr{N}}
\def\Oo{\mathscr{O}}

\def\Qq{\mathscr{Q}}

\def\Uu{\mathscr{U}}
\def\Vv{\mathscr{V}}

\def\Xx{\mathscr{X}}

\def\Ga{\Gamma}
\def\La{\Lambda}

\def\a{\alpha} 

\def\d{\delta}
\def\f{\rightarrow}

\def\h{\varphi}
\def\l{\lambda}
\def\o{\mathfrak{o}}
\def\p{\mathfrak{p}}
\def\s{\sigma}
\def\t{\theta}

\def\w{\varpi}

\def\lp{\langle}
\def\rp{\rangle}
\def\>{\geqslant}
\def\<{\leqslant}

\def\ff{\longrightarrow}

\def\ffr#1{\smash{\mathop{\ff}\limits^{#1}}}

\def\Hom{{\rm Hom}}

\def\Aut{{\rm Aut}}

\def\GL{{\rm GL}}
\def\Gal{{\rm Gal}}
\def\Ker{{\rm Ker}}
\def\Im{{\rm Im}}

\def\Res{{\rm Res}}

\def\mult#1{{#1}^{\times}}

\def\O{\Omega}
\def\Aff{{\rm Aff}}

\def\App{{\rm App}}
\def\a{a}

\def\F{k}

\def\NO{\N_{\O}}
\def\PO{\P_{\O}}
\def\UO{\U_{\O}}
\def\cpt{1}

\author{P. Delorme}
\address{Institut de Math\'ematiques de Luminy\\UMR $6206$\\
Campus de Luminy\\Case $907$\\
$13288$ Marseille Cedex $9$}
\email{delorme@iml.univ-mrs.fr}

\author{V. S\'echerre}
\address{Institut de Math\'ematiques de Luminy\\UMR $6206$\\
Campus de Luminy\\Case $907$\\
$13288$ Marseille Cedex $9$}
\email{secherre@iml.univ-mrs.fr}

\numberwithin{equation}{section} 
\theoremstyle{plain}

\title[An analogue of the Cartan decomposition]
{An analogue of the Cartan decomposition for $\boldsymbol{p}$-adic 
reductive symmetric spaces}

\begin{abstract}
Let $\F$ be a non Archimedean locally compact field of residue
characteristic different from $2$, let $\G$ be a connected reductive
group defined over $\F$, let $\s$ be an involutive $\F$-automorphism
of $\G$ and $\H$ an open $\F$-subgroup of the fixed points group of
$\s$.
We denote by $\G_\F$ (resp. $\H_\F$) the group of $\F$-points of $\G$
(resp. $\H$).
In this paper, we obtain an analogue of the Cartan decomposition for
the reductive symmetric space $\H_\F\backslash\G_\F$.
More precisely, we obtain a decomposition of $\G_\F$ as a union of
$\H_\F$-cosets which is related to the $\H_\F$-conjugacy classes of
maximal $\s$-anti-invariant $\F$-split tori in $\G$. 
When $\G$ is $\F$-split, we get a more precise result, involving the
stabilizer of a special point of the Bruhat-Tits building of $\G$ over
$\F$.
\end{abstract}

\begin{altabstract}
Soit $\F$ un corps localement compact non archim\'edien de
ca\-rac\-t\'eristique r\'esiduelle impaire, soit $\G$ un groupe
r\'eductif connexe d\'efini sur $\F$, soient $\s$ un
$\F$-automorphisme involutif de $\G$ et $\H$ un $\F$-sous-groupe
ouvert du groupe des points de $\G$ fixes par $\s$.
On note $\G_\F$ (resp. $\H_\F$) le groupe des $\F$-points de $\G$ 
(resp. $\H$).
Dans cet article, nous obtenons un analogue de la d\'ecomposition de
Cartan pour l'espace sym\'etrique r\'eductif $\H_\F\backslash\G_\F$.
Plus pr\'ecis\'ement, nous obtenons une d\'ecomposition de $\G_\F$
sous la forme d'une r\'eunion de classes modulo $\H_\F$ reli\'ee aux
classes de $\H_\F$-conjugaison de tores $\F$-d\'eploy\'es 
$\s$-anti-invariants maximaux de $\G$.
Lorsque $\G$ est d\'eploy\'e sur $\F$, nous obtenons un r\'esultat
plus pr\'ecis impliquant le stabilisateur d'un point sp\'ecial de
l'immeuble de Bruhat-tits de $\G$ sur $\F$.
\end{altabstract}

\begin{document}

\maketitle

\section*{Introduction}

Let $\F$ be a non Archimedean locally compact field of residue
characteristic different from $2$.
Let $\G$ be a connected reductive group defined over $\F$, let $\s$ be
an involutive $\F$-automorphism of $\G$ and $\H$ an open $\F$-subgroup
of the fixed points group of $\s$.
We denote by $\G_\F$ (resp. $\H_\F$) the group of $\F$-points of $\G$
(resp. $\H$).
Harmonic analysis on the reductive symmetric space
$\H_\F\backslash\G_\F$ is the study of the action of $\G_\F$ on the
space of complex square integrable functions on
$\H_\F\backslash\G_\F$.
This study is related to the classification of $\H_\F$-distinguished 
representations of $\G_\F$, that is representations having a nonzero
space of $\H_\F$-invariant linear forms.
The question of distinguishedness has been studied intensively for
$\GL_n$ and related groups. 
See for instance \cite{AP,AT,HakimMao,HM1,HM2,Prasad,Prasad2} for a
(non ex\-haustive) list of works on this question.
Some other aspects of that problem, including the Plancherel formula,
have been studied by Offen \cite{Offen} for spherical representations,
in three particular cases related to $\GL_n$.
Blanc and Delorme \cite{BD} have studied parabolically induced
representations for a general reductive symmetric space
$\H_\F\backslash\G_\F$.
In this paper, we investigate the geometry of the space 
$\H_\F\backslash\G_\F$.

Connected reductive groups over $\F$ can be considered as reductive
symmetric spaces. 
Indeed, if $\G'$ is such a group, the map $\s:(x,y)\mapsto(y,x)$
defines a $\F$-involution of the connected reductive group
$\G=\G'\times\G'$ 
whose fixed points group $\H$ is the diagonal image of $\G'$ in $\G$.
Hence the reductive symmetric space $\H_\F\backslash\G_\F$ naturally
identifies with the group $\G'_\F$. 
Moreover, if $\K'$ is a subgroup of $\G'_\F$ and if we set
$\K=\K'\times\K'$, then the $(\H_\F,\K)$-double cosets of $\G_\F$
correspond to the $\K'$-double cosets of $\G'_\F$.
In particular, if $\K'$ is the stabilizer in $\G'_\F$ of a special
point in the (reduced) Bruhat-Tits building of $\G'$ over $\F$, the
decomposition of $\H_\F\backslash\G_\F$ into $\K$-orbits corresponds
to the Cartan decomposition of $\G'_\F$ relative to $\K'$ 
(see \cite[Proposition 4.4.3]{BT}).

In this paper, we obtain an analogue of the Cartan decomposition for a
general reductive symmetric space $\H_\F\backslash\G_\F$.
In \cite{H,HH,HW} A.\ and G.\ Helminck and Wang studied two kinds of
objects which are related to our problem:
\begin{itemize}
\item[(i)]
$\H_\F$-conjugacy classes of maximal $\s$-anti-invariant
$\F$-split tori of $\G$ (called maximal {\it $(\s,\F)$-split} tori in 
\cite{H}, see also Definition \ref{Justine});
\item[(ii)]
$\H_\F$-conjugacy classes of the parabolic $\F$-subgroups $\P$ of
$\G$ which are opp\-osite to $\s(\P)$ 
(called {\it $\s$-split} parabolic $\F$-subgroups in \cite{HW} and 
{\it $\s$-parabolic} $\F$-subgroups in this paper, see Definition
\ref{Juliette}). 
\end{itemize}

Let $\{\A^{j}\ |\ j\in \J\}$ be a set of representatives of the
$\H_\F$-conjugacy classes of maximal $(\s,\F)$-split tori in $\G$.
For each $j$, we denote by $\W_{\G_{\F}}(\A^{j})$
(resp. $\W_{\H_{\F}}(\A^{j})$) the quotient of the normalizer of
$\A^{j}$ in $\G_{\F}$ (resp. in $\H_\F$) by its centralizer.
According to Helminck and Wang \cite{HW}, the set $\J$ 
is finite and, for $j\in\J$, the group $\W_{\G_{\F}}(\A^{j})$ is the
Weyl group of a root system.
Moreover, let $\A$ be a max\-imal $(\s,\F)$-split torus of $\G$, let
$\SS$ be a $\s$-stable maximal $\F$-split torus of $\G$ cont\-aining
$\A$ and $\P$ a minimal $\s$-parabolic $\F$-subgroup of $\G$ containing
$\SS$.
Then, according to \cite[Theorem 3.6]{HH}, the finite union:
\begin{equation}
\label{Samantha}
\bigcup_{j\in\J}\W_{\H_{\F}}(\A^{j})\backslash\W_{\G_{\F}}(\A^{j})
\end{equation}
classifies the open $(\H_\F,\P_\F)$-double cosets of $\G_\F$.
For each $j\in\J$, we choose:
\begin{enumerate}
\item
a set $\Nn_j\subset\N_{\G_\F}(\A^j)$ of representatives of 
$\W_{\H_{\F}}(\A^{j})\backslash\W_{\G_{\F}}(\A^{j})$;
\item
an element $y_j\in\G_{\F}$ such that $y_j\A y_j^{-1}=\A^{j}$;
\end{enumerate}
and we denote by $\Nn$ the set of all $ny_j$ for $j\in\J$ and
$n\in\Nn_j$.
Note that $\Nn$ is a set of representatives of (\ref{Samantha}).
Let $\w$ be a uniformizer of $\F$, let $\La$ be the lattice formed by
the images of $\w$ by the various algebraic one-parameter subgroups of
$\A$ and let $\La^{-}$ denote the subset of anti-dominant elements of
$\La$ relative to $\P$.
Then we can state our first main result (see Theorem
\ref{Anaximandre}):

\begin{theo}
\label{ThA}
There exists a comp\-act subset $\O$ of $\G_{\F}$ such that:
\begin{equation*}
\G_{\F}=\bigcup\limits_{n\in\Nn}\H_\F n\La^{-}\O.
\end{equation*}
\end{theo}

In order to prove this result, we make a large use of the
Bruhat-Tits theory \cite{BT,BT2}.
Let $\Bb$ be the (reduced) Bruhat-Tits building of $\G$ over $\F$.
It is endowed with an action of $\s$.
Then the proof of Theorem \ref{ThA} is based on the following result
(see Proposition \ref{Raince}):

\begin{prop}
\label{PropA}
$\Bb$ is the union of its $\s$-stable apartments.
\end{prop}

This result can be rephrased as follows.
Let $\SS$ be a $\s$-stable max\-imal $\F$-split torus of $\G$, let $\N$
its normalizer in $\G$ and let $\Oo$ be the set of all $g\in\G_\F$
such that $g^{-1}\s(g)\in\N_\F$.
Then we have $\G_\F=\Oo\K$, where $\K$ is the stabilizer in $\G_\F$ of
any point of the apartment corresp\-onding to $\SS$ 
(see Proposition \ref{Euripide}).

Let us mention that the question of the disjointness of the various 
components appearing in the decomposition of $\G_\F$ given by Theorem
\ref{ThA} has been investigated by Lagier \cite{Nat}.

When the group $\G$ is $\F$-split, we obtain a refinement of
Theorem \ref{ThA}, which is based on the following refinement of
Proposition \ref{PropA} (see Proposition \ref{Rubempre}):

\begin{prop}
\label{PropB}
Let $x$ be a special point of $\Bb$.
There is a $\s$-stable max\-imal $\F$-split torus $\SS$ of $\G$ such
that the apartment corresponding to $\SS$ cont\-ains $x$, and such
that the maximal $\s$-anti-invariant subtorus of $\SS$ is a maximal
$(\s,\F)$-split torus of $\G$.
\end{prop}

We thus obtain our second main result (see Theorem
\ref{SplendeursEtMiseres}): 

\begin{theo}
\label{ThB}
Let $\K$ be the stabilizer in $\G_\F$ of a special point in $\Bb$.
Then:
\begin{equation*}
\G_{\F}=\bigcup\limits_{j\in\J}
\H_\F y_{j}\SS_\F\K.
\end{equation*}
\end{theo}

Note that Proposition \ref{PropB} is no longer true for non-split
groups, as proven in \S\ref{CounterStrike}.

The paper is organized as follows.
In \S1 we recall the main properties of the Bruhat-Tits build\-ing
attached to a connected reductive group defined over $\F$.
In \S2 we study the set of all apartments containing a given
$\s$-stable subset of the building, and we prove Proposition
\ref{PropA}.
In \S3 we prove our first  main result (Theorem \ref{ThA}).
In \S4 we are devoted to the case where $\G$ is $\F$-split. 
We prove Proposition \ref{PropB} and Theorem \ref{ThB}.
Finally, in \S5 we study in more details the two following examples:
\begin{enumerate}
\item
$\G_{\F}=\GL_{n}(\F)$ and $\s(g)=$ transpose of $g^{-1}$.
\item
$\G_{\F}=\GL_{n}(\F')$ with $\F'$ quadratic over $\F$ and 
${\rm id}\neq\s\in\Gal(\F'/\F)$. 
\end{enumerate}
When $n=2$ and $\F'$ is totally ramified over $\F$, Example (2)
provides an example of a non-split group for which Proposition
\ref{PropB} is not satisfied.

{\it 
After this work was finished, we learnt that Y.\ Benoist and H.\ Oh 
\cite{BO} proved a result equivalent to Theorem \ref{ThA}, with a
weaker assumption on $\F$ (they only assume that its characteristic is
not $2$).
They also use the Bruhat-Tits building, but in a different way. 
}

We thanks F.\ Court\`es, B.\ Lemaire, G.\ Rousseau, S.\ Stevens for
stimulating discussions.
Particular thanks to Joseph Bernstein for having suggested to the
first author the use of the Bruhat-Tits building, and to Jean-Pierre
Labesse for having answered numerous questions, in particular about
the structure of algebraic groups. 


\section{The Bruhat-Tits building}
\label{sec1}

Let $\F$ be a non Archimedean non discrete locally compact field, and
let $\omega$ be its normalized valuation.
In this section, we recall the main properties of the
Bruhat-Tits build\-ing attached to a connected reductive group defined
over $\F$.
The reader may refer to the original construction of Bruhat-Tits
\cite{BT,BT2} or to more concise presentations \cite{Lan,SS,Tits}.  

If $\G$ is a linear algebraic group defined over $\F$, the group of
its $\F$-points will be denoted by $\G_{\F}$ or $\G(\F)$, and its
neutral component will be denoted by $\G^{\circ}$.
If $\H$ is a subset of $\G$, then $\N_{\G}(\H)$ (resp. $\Z_{\G}(\H)$) 
denotes the normalizer (resp. the centralizer) of $\H$ in $\G$.

If $\X$ is a subset of $\G$, then ${}^g\X$ denotes the left conjugate
of $\X$ by $g\in\G$. 

\subsection{}

Let $\G$ be a connected reductive group defined over $\F$, and let
$\SS$ be a maximal $\F$-split torus of $\G$.
We denote by $\X^{*}(\SS)=\Hom(\SS,\GL_1)$ (resp. by
$\X_{*}(\SS)=\Hom(\GL_1,\SS)$) the group of algebraic characters 
(resp. co\-characters) of $\SS$.
We define a map:
\begin{equation}
\label{dualiteS}
\X_{*}(\SS)\times\X^{*}(\SS)\f\ZZ
\end{equation}
as follows.
If $\l\in\X_{*}(\SS)$ and $\chi\in\X^{*}(\SS)$, 
then $\chi\circ\l$ is an endomorphism of the multiplicative group 
$\GL_{1}$, which corresponds to an endomorphism of the ring
$\ZZ[t,t^{-1}]$.
It is of the form $t\mapsto t^{n}$ for some $n\in\ZZ$.
This integer $n$ is denoted by $\lp\l,\chi\rp$.
The map (\ref{dualiteS}) defines a perfect duality 
(see \cite[\S8.6]{Borel}).

\subsection{}
\label{HoraceBianchon}

Let $\N$ (resp. $\Z$) denote the normalizer (resp. the centralizer) 
of $\SS$ in $\G$.
If we extend (\ref{dualiteS}) by $\RR$-linearity, there exists a
unique group homomorphism:
\begin{equation}
\label{nu1}
\nu:\Z_\F\f\X_{*}(\SS)\otimes_{\ZZ}\RR
\end{equation}
such that the condition:
\begin{equation}
\label{A2}
\lp\nu(z),\chi\rp=-\omega(\chi(z))
\end{equation}
holds for any $z\in\Z_\F$ and any $\F$-rational character
$\chi\in\X^{*}(\Z)_{\F}$ (see \cite[\S 1.2]{Tits}).
According to \cite[Proposition 1.2]{Lan}, the kernel of (\ref{nu1}) is
the maximal compact subgroup of $\Z_\F$.
It will be denoted by $\Z_\F^{\cpt}$. 

\begin{rema}
\label{RemCptMax}
Note that the intersection $\SS_\F\cap\Z_{\F}^{\cpt}$ is equal to the
maximal compact subgroup of $\SS_\F$, which we denote by $\SS_\F^{\cpt}$.
Indeed $\SS_\F^{\cpt}$ contains the compact subgroup
$\SS_\F\cap\Z_{\F}^{\cpt}$ of $\SS_\F$ and is contained in the maximal
compact subgroup $\Z_{\F}^{\cpt}$ of $\Z_{\F}$.
According to \cite[\S1.2]{Tits}, the quotient
$\La_{\Z}=\Z_\F/\Z_{\F}^{\cpt}$ is a free abelian group of rank
$\dim\SS$, and the image of $\SS_{\F}$ in $\La_{\Z}$ has finite
index. 
\end{rema}

\subsection{}

Let ${\rm C}$ denote the connected centre of $\G$ and let 
$\X_{*}({\rm C})$ be the group of its algebraic cocharacters.
It is a subgroup of the free abelian group $\X_{*}(\SS)$.
We denote by $\Aa$ the space: 
\begin{equation*}
\V=(\X_{*}(\SS)\otimes_{\ZZ}\RR)/(\X_{*}({\rm C})\otimes_{\ZZ}\RR)
\end{equation*}
considered as an affine space on itself and by $\Aff(\Aa)$ the group
of its affine automorphisms.
By making $\V$ act on $\Aa$ by translations, we can think to $\V$ as
a subgroup of $\Aff(\Aa)$.
It is the kernel of the natural group homomorphism
$\Aff(\Aa)\to\GL(\V)$ which associates to any affine automorphism its
linear part.

\subsection{}
\label{Hyperbate}

The map (\ref{nu1}) induces a homomorphism:
\begin{equation}
\label{nu2}
\Z_\F\f\Aff(\Aa)
\end{equation}
which is still denoted by $\nu$.
Its image is contained in $\V$.
An important property of this homomorphism is that it extends to a 
homomorphism $\N_\F\f\Aff(\Aa)$ (see \cite[\S 1.2]{Tits}). 
It does not extend in a unique way, but two homomorphisms extending
(\ref{nu2}) to $\N_\F$ are conj\-ugated by a {\it unique} element of
$\Aff(\Aa)$ (see \cite[Proposition 1.8]{Lan}).

\subsection{}
\label{consAppart}

The affine space $\Aa$ endowed with an action of $\N_\F$ defined
by a group homomorphism $\nu:\N_\F\f\Aff(\Aa)$ extending the
homomorphism (\ref{nu2}) is called the (reduced) {\it apartment}
attached to $\SS$.
It satisfies the conditions:
\begin{enumerate}
\item[{\bf A1}]
$\Aa$ is an affine space on $\V$;
\item[{\bf A2}]
$\nu$ is a group homomorphism $\N_\F\f\Aff(\Aa)$ extending the
canonical homomorphism $\Z_\F\to\V$.
\end{enumerate}
It has the following unicity property.
If $(\Aa',\nu')$ satisfy {\bf A1} and {\bf A2}, then there is a unique
affine and $\N_\F$-equivariant isomorphism from $\Aa'$ to $\Aa$.

\begin{rema}
\label{Sardanapale}
We obtain the {\it non reduced} apartment $\Aa_{\rm nr}$ by replacing 
$\V$ by $\X_{*}(\SS)\otimes_{\ZZ}\RR$.
This is the point of view of Tits \cite{Tits}.
The non reduced apartment is not as canonical as the reduced one: two
homomorphisms extending the map 
$\nu_{\rm nr}:\Z_\F\f\Aff(\Aa_{\rm nr})$ to $\N_\F$ are conjugated by
an element of $\Aff(\Aa_{\rm nr})$ which is not necessarily unique
(see \cite[\S1]{Lan} and also \cite[\S 1.2]{Tits}).
\end{rema}

\subsection{}

Let $\Phi=\Phi(\G,\SS)$ denote the set of roots of $\G$ relative to
$\SS$.
It is a subset of $\X^{*}(\SS)$.
Therefore, any root $\a\in\Phi$ can be seen as a linear form on
$\X_{*}(\SS)\otimes\RR$ which is trivial on the subspace  
$\X_{*}({\rm C})\otimes\RR$, hence as a linear form on $\V$ 
(see \cite[\S 1]{Lan}).

For $a\in\Phi$, we denote by $\U_{\a}$ the root subgroup associated to
$\a$, which is a unipotent subgroup of $\G$ normalized by $\Z$ (see
\cite[Proposition 21.9]{Borel}), and by $s_{\a}$ the reflection
corresponding to $\a$, considered as an element of $\GL(\V)$ | or,
more precisely, of the quotient of $\nu(\N_\F)$ by $\nu(\Z_\F)$.

\subsection{}
\label{SpecialeBanane}

Let $\a\in\Phi$ and $u\in\U_{\a}(\F)-\{1\}$.
The intersection: 
\begin{equation}
\label{InterMU}
\U_{-\a}(\F)u\U_{-\a}(\F)\cap\N_\F
\end{equation}
consists of a single element, called $m(u)$, whose image by $\nu$ is an
affine reflection whose linear part is $s_{\a}$ 
(see \cite[\S5]{BorelTits}).
The set $\Hh_{a,u}$ of fixed points of $\nu(m(u))$ is an affine
hyperplane of $\Aa$, which is called a {\it wall} of $\Aa$.

A {\it chamber} of $\Aa$ is a connected component of the complementary
in $\Aa$ of the union of its walls.
Note that a chamber is open in $\Aa$.

A point $x\in\Aa$ is said to be {\it special} if, for all root
$a\in\Phi$, there is a root $b\in\Phi\cap\RR_{+}a$ and an element
$u\in\U_{b}(\F)-\{1\}$ such that $x\in\Hh_{b,u}$ (see
\cite[\S1.2.3]{LanCrelle} and also \cite[\S1.9]{Tits}).

\subsection{}
\label{Chouans}

Let $\t(\a,u)$ denote the affine function $\Aa\to\RR$ whose linear
part is ${\a}$ and whose vanishing hyperplane is the wall $\Hh_{a,u}$ 
of fixed points of $\nu(m(u))$.
We fix a base point in $\Aa$, so that $\Aa$ can be identified with
the vector space $\V$.
For $r\in\RR$, we set: 
\begin{equation*}
\U_{\a}(\F)_{r}=\{u\in\U_{\a}(\F)-\{1\}\ |\
\t(\a,u)(x)\>\a(x)+r\ \text{for all}\ x\in\Aa\}\cup\{1\}. 
\end{equation*}
Thus we obtain a filtration of $\U_\a(\F)$ by subgroups.
If we change the base point in $\Aa$, this filtration is only modified
by a translation of the indexation.

\subsection{}

Let $\O$ be a nonempty subset of $\Aa$.
We set:
\begin{equation*}
\NO=\{n\in\N_\F\ |\ \nu(n)(x)=x\ \text{for all}\ x\in\O\},
\end{equation*}
and we denote by $\UO$ the subgroup of $\G_\F$ generated by all the 
$\U_{\a}(\F)_{r}$ such that the affine function $x\mapsto\a(x)+r$ is
non negative on $\O$.
According to \cite[\S12]{Lan}, this subgroup is compact in $\G_{\F}$, 
and we have $n\UO n^{-1}=\U_{\nu(n)(\O)}$ for any $n\in\N_{\F}$.
In particular $\NO$ normalizes $\UO$.

The subgroup $\PO=\NO\UO$ is open in $\G_\F$ 
(see {\it loc.cit.}, Corollary 12.12). 

\subsection{}
\label{Iwap}

Let $\Phi=\Phi^{-}\cup\Phi^{+}$ be a decomposition of $\Phi$ into
positive and negative roots.
We denote by $\U^{+}$ (resp. $\U^{-}$) the subgroup of $\G_{\F}$
generated by the $\U_{a}$ for $a\in\Phi^+$ (resp. $a\in\Phi^-$), 
and we set $\UO^+=\UO\cap\U^+$ (resp. $\UO^-=\UO\cap\U^-$).
Then the group $\PO$ has the following Iwahori decomposition:
\begin{equation}
\label{iwahoriPO}
\PO=\UO^-\UO^+\NO
\end{equation}
(see \cite[Corollary 12.6]{Lan} and also \cite[\S7.1.4]{BT}). 


\subsection{}
\label{TroisGraces}

In \cite{BT,BT2}, Bruhat and Tits associate to the apartment 
$(\Aa,\nu)$ a $\G_\F$-set $\Bb=\Bb(\G,\F)$ containing $\Aa$, called
the (reduced) {\it building} of $\G$ over $\F$ and
satisfying the following conditions:
\begin{enumerate}
\item[{\bf B1}]
The set $\Bb$ is the union of the $g\cdot\Aa$ with $g\in\G_\F$.
\item[{\bf B2}]
The subgroup $\N_{\F}$ is the stabilizer of $\Aa$ in $\G_{\F}$, and
$n\cdot x=\nu(n)(x)$ for all $x\in\Aa$ and $n\in\N_\F$.
\item[{\bf B3}]
For all $a\in\Phi$ and $r\in\RR$, the subgroup $\U_{\a}(\F)_{r}$
defined in \S\ref{Chouans} fixes the subset 
$\{x\in\Aa\ |\ a(x)+r\>0\}$ pointwise.
\end{enumerate}
The building has the following unicity property.
If $\Bb'$ is a $\G_\F$-set containing $\Aa$ and satisfying {\bf B1}, 
{\bf B2} and {\bf B3}, then there is a unique $\G_\F$-equivariant
bijection from $\Bb'$ to $\Bb$ (see \cite[\S2.1]{Tits} and also
\cite[\S1.9]{PY}). 

\subsection{}

The subsets of $\Bb$ of the form $g\cdot\Aa$ with $g\in\G_\F$ are
called {\it apartments}.
According to {\bf B1} the building is the union of its apartments.

For $g\in\G_\F$, the apartment $g\cdot\Aa$ can be naturally endowed
with a structure of affine space and an action of ${}^g\N_\F$ by
affine isomorphisms.
Upto unique iso\-morphism, it is the apartment attached to the maximal
$\F$-split torus ${}^g\SS$ (see \S\ref{consAppart}).
This defines a unique $\G_\F$-equivariant map:
\begin{equation}
\label{Herrera}
\G\supset\SS'\mapsto\Aa(\SS')\subset\Bb
\end{equation}
between maximal $\F$-split tori of $\G$ and apartments of $\Bb$, such
that $\SS$ maps to $\Aa$.

Note that $\Bb$ does not depend on the maximal $\F$-split torus
$\SS$.
Indeed, let $\SS'$ be a maximal $\F$-split torus of $\G$, let
$(\Aa',\nu')$ be the apartment attached to $\SS'$ and $\Bb'$ 
the building of $\G$ over $\F$ relative to this apartment (see 
\S\ref{TroisGraces}).
If we identify $\Aa'$ with the unique apartment of $\Bb$
corresponding to $\SS'$ {\it via} (\ref{Herrera}), then $\Bb'=\Bb$.

\subsection{}
\label{TroisParques}

The building has the following important properties
(see \cite[\S7.4]{BT} and also \cite[\S13]{Lan}): 
\begin{enumerate}
\item
Let $\O$ be a nonempty subset of $\Aa$.
Then $\PO$ is the subgroup of $\G_{\F}$ made of all elements
fixing $\O$ pointwise.
\item
Let $g\in\G_{\F}$.
There exists $n\in\N_{\F}$ such that $g\cdot x=n\cdot x$ for any
element $x\in\Aa\cap g^{-1}\cdot\Aa$.
\end{enumerate}
In particular, Property (1) together with {\bf B2} imply that
$\NO=\N_{\F}\cap\PO$.

\subsection{}
\label{ActionSigma}

Let $\s$ be a $\F$-automorphism of $\G$.
There is a unique bijective map from $\Bb$ to itself, which we still
denote by $\s$, such that:
\begin{itemize}
\item[(i)]
the condition:
\begin{equation*}
\label{act}
\s(g\cdot x)=\s(g)\cdot\s(x)
\end{equation*}
holds for any $g\in\G_{\F}$ and $x\in\Bb$;
\item[(ii)]
the map $\s$ permutes the apartments and, for any apartment $\Aa$, 
the restriction of $\s$ to $\Aa$ is an affine isomorphism from $\Aa$ 
onto its image.
\end{itemize}
This gives us an action of the group $\Aut_{\F}(\G)$ of
$\F$-automorphisms of $\G$ on the building (see \cite[\S 4.2.12]{BT2}).

\subsection{}
\label{TheExtendedBuilding}

Let $\V^1$ denote the dual space of $\X^*(\G)\otimes\RR$. 
The {\it extended building} of $\G$ over $\F$ is the product 
$\Bb^1=\Bb\times\V^1$, where $\G_\F$ acts on $\V^1$ by: 
\begin{equation*}
\label{ActionBExt}
g\cdot\chi=-\omega(\chi(g)),
\end{equation*}
for any $g\in\G_\F$ and any $\F$-rational character
$\chi\in\X^*(\G)_{\F}$.
The $\G_\F$-stabilizer of the reduced building $\Bb\times\{0\}$, 
considered as a subset of the extended build\-ing $\Bb^1$, is denoted
by $\G_{\F}^{\cpt}$.
It is the subgroup of all $g\in\G_\F$ such that $\omega(\chi(g))=0$
for any $\chi\in\X^*(\G)_{\F}$.


\begin{rema}
\label{ReseauxGD}
Let $\D$ denote the maximal $\F$-split torus of the connected centre
$\C$ of $\G$.
Then the quotient $\La_{\G}=\G_\F/\G^{\cpt}_\F$ is a free abelian
group of rank $\dim\D$, and the image of $\D_{\F}$ in $\La_{\G}$ has
finite index. 
\end{rema}

The action of $\Aut_{\F}(\G)$ on $\G$ induces an action of
$\Aut_{\F}(\G)$ on $\V^1$, hence on the extended building $\Bb^1$.

\begin{rema}
Let $\Ga$ be a finite subgroup of $\Aut_{\F}(\G)$ whose order is prime
to the residue characteristic of $\F$, and let $\H$ be the neutral
component of the fixed points subgroup $\G^{\Ga}$.
Prasad and Yu \cite{PY} proved the existence of a $\H_{\F}$-equivariant
map $\iota:\Bb^{1}(\H,\F)\to\Bb^{1}(\G,\F)$ whose image is the set
of $\Ga$-invariant points.
Moreover, such a map is {\it toral} in the sense of \cite{Lan}, 
which means that for any maximal $\F$-split torus $\T$ of $\H$, there
is a maximal $\F$-split torus $\SS$ of $\G$ containing $\T$ such that 
$\iota$ maps $\Aa_{\rm nr}(\H,\T)$ to $\Aa_{\rm nr}(\G,\SS)$ by an
affine transformation (see \cite[Theorem 1.9]{PY}, and Remark
\ref{Sardanapale} for the definition of the non reduced apartment
$\Aa_{\rm nr}$).
\end{rema}


\section{Existence of $\s$-stable apartments}
\label{sec3}

From now on, $\F$ will be a non-Archimedean locally compact
field of residue characteristic different from $2$.
Let $\G$ be connected reductive group defined over $\F$ and let $\s$
be a $\F$-involution on $\G$.

According to \S\ref{ActionSigma}, the building $\Bb$ of $\G$
over $\F$ is endowed with an action of $\s$.
In this section, we prove that, for any $x\in\Bb$, there exists a
$\s$-stable apartment containing $x$.
We keep using notations of Section \ref{sec1}.

\subsection{}
\label{Enmonameetconscience}

Let $\O$ be a nonempty $\s$-stable subset of $\Bb$ contained in some
apartment, and let $\App(\O)$ be the set of all apartments of $\Bb$
containing $\O$. 
It is a nonempty set on which the group $\PO$ acts transitively 
(see \cite[Corollary 13.7]{Lan}).
Because $\O$ is $\s$-stable, both $\PO$ and $\App(\O)$ are
$\s$-stable.
Note that the $\s$-stable apartments containing $\O$ are exactly the
$\s$-invariant points in $\App(\O)$.

\subsection{}
\label{toune}

Let us fix an apartment $\Aa\in\App(\O)$ and an element $u\in\PO$
such that $\s(\Aa)=u\cdot\Aa$.
Let $\N$ denote the normalizer in $\G$ of the maximal $\F$-split torus
of $\G$ corresponding to $\Aa$.
As $\s$ is involutive, we have:
\begin{equation}
\label{condg}
\s(u)u\in\PO\cap\N_{\F}=\NO.
\end{equation}
The map $\rho:g\mapsto g\cdot\Aa$ induces a $\PO$-equivariant
bijection between the homogeneous spaces $\PO/\NO$ and $\App(\O)$.
The automorphism:
\begin{equation}
\t:x\mapsto u^{-1}\s(x)u
\end{equation}
of the group $\G_{\F}$ stabilises $\PO$ and $\NO$.
Indeed $\s(\N_{\F})=u\N_{\F}u^{-1}$, and:
\begin{equation*}
\t(\NO)=u^{-1}\s(\PO\cap\N_{\F})u=\PO\cap u^{-1}\s(\N_{\F})u=\NO.
\end{equation*}
Note that the condition (\ref{condg}) means that $\t\circ\t$ is
conjugation by some element of $\NO$.
As $\NO$ is $\t$-stable, the map:
\begin{equation}
\label{ActionPN}
(\s,g\NO)\mapsto u\t(g\NO),\quad g\in\PO,
\end{equation}
defines an action of $\s$ on $\PO/\NO$, making $\rho$ into a
$\s$-equivariant bijection.
Note that this action differs from the natural action of $\s$ on
$\PO/\NO$ (which obviously has fixed points). 

\subsection{}

Let $\O$ be a nonempty $\s$-stable subset of $\Bb$ contained in some
apartment. 

\begin{prop}
\label{AppartStableNoWall}
Assume that $\O$ contains a point of a chamber of $\Bb$.
Then $\O$ is contained in some $\s$-stable apartment.
\end{prop}

\begin{proof}
First we describe $\PO/\NO$ as a projective limit of finite
$\s$-sets.
According to \cite[\S1.2]{Cartier}, Example (f), the group $\G_{\F}$
is locally compact and totally disconnected.
Therefore we can choose a decreasing filtration $(\Q^{i})_{i\>0}$ of
the open subgroup $\PO$ of $\G_{\F}$ satisfying the following
properties: 
\begin{enumerate}
\item[(A)]
The intersection of the $\Q^{i}$ is reduced to $\{1\}$.
\item[(B)]
For any $i\>0$, the subgroup $\Q^{i}$ is compact open and normal in
$\PO$.
\end{enumerate}
For $i\>0$, let $\PO^i$ denote the intersection of the subgroups
$\NO\Q^{i}$ and $\t(\NO\Q^{i})$. 
The $\PO^{i}$ form a decreasing filtration of $\PO$, and we claim that
such a filtration satisfies the following properties: 
\begin{enumerate}
\item
The intersection of the $\PO^{i}$ is reduced to $\NO$.
\item
For any $i\>0$, the subgroup $\PO^{i}$ is $\t$-stable
and of finite index in $\PO$.
\end{enumerate}
As $\NO$ is $\t$-stable, it is contained in the intersection of the
$\PO^{i}$.
Let $g$ be in this intersection.
For any $i\>0$, there exist $n_{i}\in\NO$ and
$q_{i}\in\Q^{i}$ such that $g=n_{i}q_{i}$.
Because of Property (A) above, $q_{i}$ converges to $1$.
Therefore $n_{i}$ converges to a limit contained in the closed subgroup
$\NO$, and this limit is $g$.
This proves Property (1).

Now recall that $\t\circ\t$ is conjugation by some element of $\NO$. 
This implies that $\PO^{i}$ is $\t$-stable.
As $\PO^{i}$ is open in $\PO$ and contains $\NO$, the quotient
$\PO/\PO^{i}$ can be identified with the quotient of $\UO$, which is
compact, by some open subgroup.
This gives the expected result.

Because of Property $(2)$, the map:
\begin{equation*}
(\s,g\PO^{i})\mapsto u\t(g\PO^{i}),\quad g\in\PO,
\end{equation*}
defines an action of $\s$ on the finite quotient
$\PO/\PO^{i}$.
We get a projective system $(\PO/\PO^{i})_{i\>0}$ of finite 
$\s$-sets.
Because $\PO$ is complete, and thanks to Property (1), the natural
$\s$-equivariant map from $\PO/\NO$ to the projective limit
of the $\PO/\PO^{i}$ is bijective.

\begin{lemm}
\label{LimPro}
Let $(\Xx^{i})_{i\>0}$ be a projective system of finite $\s$-sets, and
let $\Xx$ be its projective limit. 
Assume that, for each $i\>0$, the cardinal of $\Xx^i$ is odd.
Then $\Xx$ has a $\s$-invariant point. 
\end{lemm}

\begin{proof}
Because each of the $\Xx^{i}$ has an odd cardinal, each of them
contains a $\s$-invariant element.
Suppose that we have constructed for some $i\>1$ a $\s$-invariant
element $x_i\in\Xx^i$.
The fiber of $x_i$ in $\Xx^{i+1}$ is $\s$-stable and its cardinal is
the quotient of the cardinal of $\Xx^{i+1}$ by the one of $\Xx^{i}$.
Therefore it is odd.
We deduce from this that there exists a $\s$-invariant element
$x_{i+1}\in\Xx^{i+1}$ whose image in $\Xx^{i}$ is $x_i$.
By induction, we get a $\s$-invariant element $x\in\Xx$.
\end{proof}

Let $p$ denote the residue characteristic of $\F$.
Recall that $p$ is assumed to be odd.

\begin{lemm}
\label{lemmppower}
Let $\K$ be a normal subgroup of finite index in $\PO$ containing
$\NO$. 
Then the index of $\K$ in $\PO$ is a power of $p$. 
\end{lemm}

\begin{proof}
Let $\SS$ be the maximal $\F$-split torus associated to $\Aa$, let 
$\Phi$ denote the set of roots of $\G$ relative to $\SS$ and let
$\Phi=\Phi^{-}\cup\Phi^{+}$ be a decomposition of $\Phi$ into positive
and negative roots.
According to \S\ref{Iwap}, the group $\PO$ has the following Iwahori
decomposition:
\begin{equation}
\PO=\UO^-\UO^+\NO.
\end{equation}
The fact that $\O$ contains a point of a chamber of $\Bb$ implies 
that the group $\NO$ is reduced to $\Ker(\nu)$,
hence normalizes the groups $\UO^+$ and $\UO^-$.
The index of $\K$ in $\PO$ can be decomposed as follows:
\begin{equation}
\label{Pingouin}
(\PO:\K)=(\PO:\UO^+\K)\cdot(\UO^+\K:\K).
\end{equation}
In a first hand, the index $(\UO^+\K:\K)=(\UO^+:\UO^+\cap\K)$ is a
power of $p$, because $\UO^+$ is a pro-$p$-group ({\it i.e.} a
projective limit of finite discrete $p$-groups).
In the other hand, the index $(\PO:\UO^+\K)$ is equal to 
$(\UO^-:\UO^-\cap\UO^+\K)$, which is a power of $p$ because
$\UO^-$ is a pro-$p$-group.
The result follows.
\end{proof}

According to Lemma \ref{lemmppower}, the cardinal of each set
$\PO/\PO^{i}$ with $i\>0$ is odd.
Proposition \ref{AppartStableNoWall} now follows from Lemma
\ref{LimPro}.
\end{proof}

\subsection{}
\label{BrandeburgerTor}

We now prove the main result of this section.

\begin{prop}
\label{Raince}
For any $x\in\Bb$, there exists a $\s$-stable apartment containing
$x$. 
\end{prop}

\begin{proof}
Let $x$ be a point in $\Bb$, and let $y$ be a point of a chamber of
$\Bb$ whose adherence contains $x$.
The set $\O=\{y,\s(y)\}$ is a $\s$-stable subset of $\Bb$ satisfying
the conditions of Proposition \ref{AppartStableNoWall}.
Hence we get a $\s$-stable apartment of $\Bb$ containing $y$.
Such an apartment contains the adherence of the chamber of $y$.
In particular, it contains $x$.
\end{proof}


\section{Decomposition of $\H_{\F}\backslash\G_{\F}$}
\label{AubergeDuPilori}

Let $\F$ be a non-Archimedean locally compact field of residue
characteristic different from $2$. 
Let $\G$ be a connected reductive group defined over $\F$, let
$\s$ be an involutive $\F$-automorphism of $\G$ and let $\H$ be an
open $\F$-subgroup of the fixed points group $\G^{\s}$. 
Equivalently, $\H$ is a $\F$-subgroup of $\G^{\s}$ containing the
neutral component $(\G^{\s})^{\circ}$ (see \cite{BDS}).

\subsection{}

Let $\SS$ be a maximal $\F$-split torus of $\G$, and let $\Aa$ denote
the corresponding apartment.

\begin{lemm}
\label{AppTore}
$\Aa$ is $\s$-stable if, and only if $\SS$ is $\s$-stable.
\end{lemm}

\begin{proof}
This comes from the fact that the apartment corresponding to
$\s(\SS)$ is the image of $\Aa$ by $\s$. 
\end{proof}

\subsection{}
\label{IciOEstDefini}

We now assume that $\SS$ is $\s$-stable.
Let $\N$ (resp. $\Z$) denote the normalizer (resp. the centralizer) 
of $\SS$ in $\G$.
Let $\Oo=\Oo_\SS$ denote the set of all $g\in\G_\F$ such that
$g^{-1}\s(g)\in\N_{\F}$.

\begin{prop}
\label{prop34}
$\Oo$ is a finite union of $(\H_{\F},\Z_{\F})$-double cosets.
\end{prop}

\begin{proof}
Let us fix a minimal parabolic $\F$-subgroup $\P$ of $\G$ containing
$\SS$.
According to \cite[Proposition 6.8]{HW}, the map
$g\mapsto\H_{\F}g\P_{\F}$ induces a bijection between the
$(\H_{\F},\Z_{\F})$-double cosets in $\Oo$ and the
$(\H_{\F},\P_{\F})$-double cosets in $\G_{\F}$.
The result then follows from \cite[Corollary 6.16]{HW}.

We now give a direct proof of this result.
We have an exact sequence:
\begin{equation*}
\G_{\F}^{\s}=\H^{0}(\G_{\F})\f\H^0(\G_{\F}/\N_{\F})
\ffr{\d}\H^1(\N_{\F})\f\H^1(\G_{\F}),
\end{equation*}
where $\H^0$ and $\H^1$ denote respectively the set of $\s$-fixed
points and the first set of nonabelian cohomology of $\s$ (see
\cite[Chapter I, \S5]{SerreCG}).
The transition map $\d$ induces an injective map from
$\G_{\F}^{\s}\backslash\H^0(\G_{\F}/\N_{\F})$, which is the set of
$(\G_{\F}^{\s},\N_{\F})$-double cosets of $\Oo$, into $\H^1(\N_{\F})$. 
Because $\Z_{\F}$ (resp. $\H_{\F}$) is of finite index in $\N_{\F}$
(resp. in $\G_{\F}^{\s}$), the finiteness of the number of 
$(\G_{\F}^{\s},\N_{\F})$-double cosets of $\Oo$ is equivalent to the
finiteness of the number of $(\H_{\F},\Z_{\F})$-double cosets of $\Oo$.
Therefore, it will be enough to prove that $\H^1(\N_{\F})$ is finite. 

Let $\M$ be a group with an action of $\s\in\Aut(\M)$, and let
$\M'$ be a $\s$-stable normal subgroup of $\M$.
We can form the following exact sequence:
\begin{equation*}
\H^1(\M')\f\H^1(\M)\f\H^{1}(\M/\M'),
\end{equation*}
which proves that the finiteness of $\H^1(\M')$ and $\H^{1}(\M/\M')$ 
implies the finiteness of $\H^1(\M)$.
Therefore we are reduced to proving that:
\begin{equation*}
\H^1(\N_{\F}/\Z_{\F}),\quad\H^1(\Z_{\F}/\Z_{\F}^{\cpt}),
\quad\H^1(\Z_{\F}^{\cpt}) 
\end{equation*}
are finite sets.
Recall (see \S\ref{HoraceBianchon}) that $\Z_{\F}^{\cpt}$ denotes the
maximal compact subgroup of $\Z_\F$.
Because $\N_{\F}/\Z_{\F}$ is finite, the first case is immediate.
Next, the quotient $\La=\Z_{\F}/\Z_{\F}^{\cpt}$ is a finitely
generated free abelian group. 
We have an exact sequence:
\begin{equation*}
\H^1(2\La)\ffr{\a}\H^1(\La)\f\H^1(\La/2\La).
\end{equation*}
Let $2m\in2\La$ be a cocycle, that is $2m+\s(2m)=0$, and consider it
as a cocycle in $\La$.
The identity $2m=m-\s(m)$ implies that the class $\a(2m)$ is trivial
in $\H^1(\La)$, hence that the map $\a$ is null.
Therefore $\H^1(\La)$ is embedded in $\H^1(\La/2\La)$, which is finite
because $\La/2\La$ is finite.

Now we treat the case of the compact subgroup $\Z_{\F}^{\cpt}$.
Let $\M$ be an open pro-$p$-subgroup of $\Z_{\F}^{\cpt}$.
(Its existence is a topological property of $\G_{\F}$ asserted in
\cite[\S1.2]{Cartier}, Example (f).) 
The normalizer of $\M$ in $\Z_{\F}^{\cpt}$ is open, hence of finite
index, in $\Z_{\F}^{\cpt}$.
We can therefore assume that $\M$ is normal (if not, we replace it by
the intersection of the finitely many ${}^{g}\M$ with
$g\in\Z_{\F}^{\cpt}$). 
Moreover, we assume that $\M$ is stable by $\s$ (if not, we replace it
by $\M\cap\s(\M)$). 
Then $\H^1(\M)$ is trivial because $\M$ is a pro-$p$-group and $p$ is
odd, and $\H^1(\Z_{\F}^{\cpt}/\M)$ is finite because $\M$ is of finite
index in $\Z_{\F}^{\cpt}$. 
The finiteness of $\H^1(\Z_{\F}^{\cpt})$ follows.
This ends our alternative proof of Lemma \ref{prop34}.
\end{proof}

\subsection{}

Let $\Aa$ denote the $\s$-stable apartment corresponding to $\SS$. 

\begin{lemm}
\label{gO}
We have $g\in\Oo$ if and only if $g\cdot\Aa$ is $\s$-stable.
\end{lemm}

\begin{proof}
As $\Aa$ is $\s$-stable, the apartment $g\cdot\Aa$ is
$\s$-stable if and only if $\s(g)\cdot\Aa=g\cdot\Aa$.
This amounts to saying that $g^{-1}\s(g)\in\N_{\F}$.
\end{proof}

For $x\in\Aa$, let $\P_x$ denote the subgroup $\P_{\O}$ 
(see \S\ref{Iwap}) with $\O=\{x\}$.

\begin{prop}
\label{Euripide}
Let $x$ be in $\Aa$.
Then we have $\G_\F=\Oo\P_{x}$.
\end{prop}

\begin{proof}
For $g\in\G_\F$, we set $x'=g\cdot x$. 
According to Proposition \ref{Raince}, there is a $\s$-stable
apartment $\Aa'$ containing $x'$.
Let $g'\in\Oo$ be such that $\Aa'=g'\cdot\Aa$.
According to Property (2) of \S\ref{TroisParques}, there exists
$n\in\N_{\F}$ such that we have $g'^{-1}g\cdot x=n\cdot x$. 
Hence we get $g\in\Oo\N_{\F}\P_{x}$.
As $\Oo\N_{\F}=\Oo$, we obtain the expected result.
\end{proof}

\subsection{}
\label{AlmostTrue}

If $\T$ is a $\s$-stable torus in $\G$, we denote by $\T^+$
(resp. $\T^-$) the neutral comp\-onent of $\T\cap\H$ (resp. of the
subgroup $\{t\in\T\ |\ \s(t)=t^{-1}\}$).
Note that, as $\T^+$ is open in the fixed points subgroup $\T^{\s}$,
we have $\T^+=(\T^{\s})^{\circ}$.
The torus $\T$ is the almost direct product (see \cite[xi]{Borel}) of
$\T^+$ and $\T^-$, which means that $\T$ is equal to the product
$\T^+\T^-$ and the intersection $\T^+\cap\T^-$ is finite.

\begin{defi}[Helminck-Wang \cite{HW}, \S4.4]
\label{Justine}
A $\s$-stable torus $\T$ of $\G$ is said to be {\it $(\s,\F)$-split}
if it is $\F$-split and if $\T=\T^-$.
\end{defi}

Let us recall (see \cite[Proposition 10.3]{HW}) that two arbitrary
maximal $(\s,\F)$-split tori of $\G$ are $\G_\F$-conjugated.

\subsection{}
\label{Poil}

Let $\T$ be a $\F$-split torus of $\G$, and let $\T_\F^{\cpt}$ denote
its maximal compact subgroup.
Let $\w$ be a uniformizer of $\F$.
The images of $\w$ by the various alg\-ebraic cocharacters of $\T$
form a $\s$-stable lattice in $\T_{\F}$, which will be denoted by
$\La(\T_{\F})$.

\begin{lemm}
\label{Borely}
\begin{itemize}
\item[(i)]
$\T_{\F}$ is the direct product of $\La(\T_{\F})$ and $\T_\F^{\cpt}$.
\item[(ii)]
For any $g\in\G_\F$, we have $\La({}^g\T_\F)={}^g\La(\T_{\F})$. 
\item[(iii)]
The subgroup generated by $\La(\T^+_\F)$ and $\La(\T^-_\F)$ has 
finite index in $\La(\T_{\F})$.
\end{itemize}
\end{lemm}

\begin{proof}
Only (iii) is not immediate.
First note that, as $\F$ is a non Archi\-me\-dean locally comp\-act
field of characteristic different from $2$, the subgroup of squares of
$\mult\F$ is of finite index in $\mult\F$.
This implies that $\T_\F^2=\{t^2\ |\ t\in\T_\F\}$ is of finite index
in $\T_\F$.

For any $t\in\T_\F$, the element $t^2$ can be decomposed as the
product of $t\s(t)\in\T^+_\F$ and $t\s(t)^{-1}\in\T^-_\F$.
Indeed the image of $\T$ by the map $t\mapsto t\s(t)$ is connected and
contained in $\T^{\s}$, thus in $\T^{+}$.
By a similar argument, the image of $\T$ by $t\mapsto t\s(t)^{-1}$
is contained in $\T^{-}$.

Therefore $\T_\F^2$ is contained in $\T^+_\F\T^-_\F$, thus
there is some finite subset $\Ff$ of $\T_\F$ such that
$\T_\F=\T^+_\F\T^-_\F\Ff$.
According to (i), this gives:
\begin{eqnarray*}
\label{Walalal}
\La(\T_\F)\T_\F^{\cpt}
&=&\La(\T_\F^+)(\T_\F^+)^{\cpt}\La(\T_\F^-)(\T_\F^-)^{\cpt}\Ff\\
&=&\La(\T_\F^+)\La(\T_\F^-)\T_\F^{\cpt}\Ff.
\end{eqnarray*}
We obtain the expected result by computing the quotient of this
equality by the subgroup $\La(\T_\F^+)\La(\T_\F^-)\T_\F^{\cpt}$.
\end{proof}

\subsection{}

Let $\{\A^{j}\ |\ j\in\J\}$ be a set of representatives of the
$\H_\F$-conjugacy classes of maximal $(\s,\F)$-split tori in $\G$.
We denote by $\W_{\G_{\F}}(\A^{j})$ (resp. $\W_{\H_{\F}}(\A^{j})$) the
quotient of the normalizer of $\A^{j}$ in $\G_{\F}$ (resp. in $\H_\F$)
by its centralizer.
According to \cite[Proposition 5.9]{HW}, the group
$\W_{\G_{\F}}(\A^{j})$ is the Weyl group of a root system.
In particular, it is a finite group.
(If $\s$ is trivial on the isotropic factor of $\G$ over $\F$, then
this group is trivial.)

\begin{defi}
\label{Juliette}
A parabolic subgroup $\P$ of $\G$ is said to be $\s$-{\it parabolic}
if it is opposite to $\s(\P)$, that is if $\P\cap\s(\P)$ is a Levi
subgroup of $\P$ and $\s(\P)$.
\end{defi}

\begin{rema}
This differs from the terminology used in \cite{HW}, where such
parabolic subgroups are said to be $\s$-{\it split}.
\end{rema}

\subsection{}

Let $\A$ be a max\-imal $(\s,\F)$-split torus of $\G$.

\begin{lemm}
There is a $\s$-stable maximal $\F$-split torus of $\G$ containing
$\A$.
\end{lemm}

\begin{proof}
Let $\G'$ denote the neutral component of the centralizer of $\A$ in
$\G$.
It is a connected reductive $\F$-group.
Let $\SS$ be a $\s$-stable maximal $\F$-split torus of $\G'$, whose
existence is asserted by Proposition \ref{Raince} and Lemma
\ref{AppTore} together.
Such a torus $\SS$ is a $\s$-stable maximal $\F$-split torus of $\G$ 
containing $\A$.
\end{proof}

Let $\SS$ be a $\s$-stable maximal $\F$-split torus of $\G$ containing
$\A$ and $\P$ a minimal $\s$-parabolic $\F$-subgroup of $\G$
containing $\SS$ (see \cite[\S4]{HW}).
Let $\w$ be a uniformizer of $\F$, set $\La=\La(\A_\F)$ and let
$\La^{-}$ denote the subset of anti-dominant elements of $\La$
relative to $\P$.

\begin{theo}
\label{Anaximandre}
For $j\in\J$, let $\Nn_j\subset\N_{\G_\F}(\A^j)$ be a set of 
representatives of 
$\W_{\H_{\F}}(\A^{j})\backslash\W_{\G_{\F}}(\A^{j})$
and $y_j\in\G_{\F}$ such that ${}^{y_{j}}\A=\A^{j}$.
There exists a comp\-act subset $\O$ of $\G_{\F}$ such that:
\begin{equation*}
\G_{\F}=\bigcup\limits_{j\in\J}\bigcup\limits_{n\in\Nn_j}
\H_\F ny_{j}\La^{-}\O.
\end{equation*}
\end{theo}

\begin{proof}
First let $\{u_i\ |\ i\in\I\}$ be a set of representatives of
$(\H_\F,\Z_{\F})$-double cosets in $\Oo$.
According to Lemma \ref{prop34}, such a set is finite.
Let $\Aa$ denote the apartment corresponding to $\SS$, and let
$\K$ be the stabilizer of $x$ in $\G_\F$.
Then Proposition \ref{Euripide} can be rephrased as follows: 
\begin{equation}
\label{Lorenza}
\G_{\F}=\bigcup\limits_{i\in\I}\H_{\F}u_i\Z_{\F}\K.
\end{equation}
Let $\K^{\cpt}$ denote the intersection $\K\cap\G_\F^{\cpt}$
(see \S\ref{TheExtendedBuilding}).
It is the maximal compact sub\-group of $\K$.

\begin{lemm}
\label{Lorenzzatio}
We have $\Z_{\F}\K=\La(\SS_{\F})\Ff\K^{\cpt}$ for some finite subset
$\Ff\subset\G_{\F}$. 
\end{lemm}

\begin{proof}
First note that $\Z_{\F}\cap\K^{\cpt}=\Z_{\F}^{\cpt}$.
Indeed, any element of the group $\Z_{\F}^{\cpt}$, which is the kernel
of (\ref{nu1}), acts trivially on $x$.
Therefore $\Z_{\F}^{\cpt}$ is contained in $\K$, hence in its maximal
compact sub\-group $\K^{\cpt}$.
Inversely, the compact group $\Z_{\F}\cap\K^{\cpt}$ is contained in
$\Z_\F$, hence in its maximal compact sub\-group $\Z_{\F}^{\cpt}$.
According to Remark \ref{RemCptMax}, the group 
$\SS_{\F}\Z_{\F}^{\cpt}$ has finite index in $\Z_{\F}$.
Thus there exists a finite subset $\Ff_1\subset\Z_{\F}$ 
such that $\Z_{\F}=\Ff_1\SS_{\F}(\Z_{\F}\cap\K^{\cpt})$.

Let $\D$ denote the maximal $\F$-split torus of the connected centre
$\C$ of $\G$.
According to Remark \ref{ReseauxGD}, the image of $\D_{\F}$ in
$\G_\F/\G^{1}_{\F}$ has finite index, thus its image in $\K/\K^{1}$
too.
This implies that $\D_{\F}\K^{\cpt}=\La(\D_{\F})\K^{\cpt}$ has finite
index in $\K$, thus that there exists a finite subset $\Ff_2\subset\K$
such that $\K=\Ff_2\La(\D_{\F})\K^{\cpt}$. 

Finally, we have:
\begin{eqnarray*}
\Z_{\F}\K&=&\Ff_1\SS_{\F}\K\\
&=&\Ff_1\La(\SS_{\F})\K\\
&=&\Ff_1\La(\SS_{\F})\Ff_2\La(\D_{\F})\K^{\cpt}
\end{eqnarray*}
which gives the expected result with $\Ff=\Ff_1\Ff_2$.
\end{proof}

For $i\in\I$, we set $\SS^i={}^{u_{i}}\SS$.
According to Lemmas \ref{Lorenzzatio} and \ref{Borely}(iii), there are
finite subsets $\Ff\subset\G_\F$ and $\Vv_i\subset\La(\SS^{i}_{\F})$,
for $i\in\I$, such that:
\begin{equation}
\label{Bargeton}
\H_{\F}u_i\Z_{\F}\K=\H_{\F}\La(\SS^{i-}_{\F})\Vv_iu_i\Ff\K^{\cpt}.
\end{equation}
According to \cite[Lemma 2.2]{HH}, the $(\s,\F)$-split torus
$\SS^{i-}$ is $\H_\F$-conjugated to a subtorus of $\A^j$ for some
$j\in\J$. 
We can therefore assume that, for a suitable choice of the
representative $u_{i}$, the $(\s,\F)$-split torus $\SS^{i-}$ is
contained in $\A^j$ for some $j\in\J$.
For $j\in\J$, let $\Uu_j$ be the union of the $\Vv_iu_i\Ff$ such
that $\A^{j}$ contains $\SS^{i-}$.
Together with (\ref{Lorenza}) and (\ref{Bargeton}), this gives: 
\begin{equation}
\label{LouisLambert}
\G_{\F}=
\bigcup\limits_{j\in\J}\H_{\F}\La(\A^{j}_{\F})\Uu_j\K^{\cpt}.
\end{equation}
For $j\in\J$, we fix a set $\Nn_{\H_{\F},j}$ of representatives of 
$\W_{\H_{\F}}(\A^{j})$ and we denote by $\tilde\Nn_{j}$ the set 
$\{hn\ |\ h\in\Nn_{\H_{\F},j}, n\in\Nn_{j}\}$.
It is a set of representatives of $\W_{\G_{\F}}(\A^{j})$.
From (\ref{LouisLambert}) we have:
\begin{equation}
\label{Lousteau}
\G_{\F}=\bigcup\limits_{j\in\J}\bigcup\limits_{n\in\tilde\Nn_{j}}
\H_{\F}n\La(\A^{j}_{\F})^-n^{-1}\Uu_j\K^{\cpt},
\end{equation}
where $\La(\A^{j}_{\F})^-$ denotes the subset of anti-dominant
elements of $\La(\A^{j}_{\F})$ relative to the parabolic subgroup
${}^{y_{j}}\P$.
If we remark that $\La(\A^{j}_{\F})^-={}^{y_{j}}\La^{-}$, and if we
denote by $\O$ the union of the $y_j^{-1}n^{-1}\Uu_j\K^{\cpt}$
for $j\in\J$ and $n\in\tilde\Nn_{j}$, then (\ref{Lousteau}) becomes:
\begin{equation}
\G_{\F}=\bigcup\limits_{j\in\J}\bigcup\limits_{n\in\Nn_{j}}
\H_{\F}ny_{j}\La^-\O.
\end{equation}
This gives us the expected result.
\end{proof}


\section{The split case}

In this section, we keep using notations of Section
\ref{AubergeDuPilori}.
Moreover, we assume that the reductive group $\G$ is split over $\F$.
Therefore, for any root $\a$ of $\G$ relative to some maximal
$\F$-split torus of $\G$, the root subgroup $\U_{\a}$ is
$\F$-isomorphic to the additive group.

The main results of this section are Proposition \ref{Rubempre} and
Theorem \ref{SplendeursEtMiseres}. 

\subsection{}
\label{Training}

Let $\SS$ be a $\s$-stable maximal $\F$-split torus of $\G$, let $\Aa$
be the apartment corresponding to $\SS$ and $\Phi$ the set of roots of
$\G$ relative to $\SS$.

Let $x\in\Aa$ be a special point, and let $\U_{x}$ denote the subgroup
$\U_{\O}$ (see \S\ref{Iwap}) with $\O=\{x\}$.
Let $\a\in\Phi$ be a $\s$-invariant root, which means that
$\a\circ\s=\a$.

\begin{lemm}
\label{Almani}
Assume that $\U_{-a}(\F)$ is contained in 
$\{g\in\G_\F\ |\ \s(g)=g^{-1}\}$.
Then there are $n\in\N_\F$ and $c\in\U_x$ such that $n=c^{-1}\s(c)$
and $\nu(n)$ is the affine reflection of $\Aa$ which let $x$ invariant
and whose linear part is $s_a$.
\end{lemm}

\begin{proof}
We fix a base point in the apartment $\Aa$, so that it can be
id\-ent\-ified with the vector space $\V$.
For any $b\in\Phi$, this defines a filtration of the group $\U_{b}(\F)$
(see \S\ref{Chouans}).
For $u\in\U_{b}(\F)-\{1\}$, we denote by $\h_{b}(u)$ the greatest real
number $r\in\RR$ such that $u\in\U_{b}(\F)_{r}$.

Let us choose $w\in\U_{-a}(\F)-\{1\}$ such that $x$ is contained in
the wall $\Hh_{-a,w}$. 
Thus $\nu(m(w))$ is the affine reflection of $\Aa$ which fixes $x$
and whose linear part is $s_a$, and we can set:
\begin{equation*}
n=m(w)\in\N_\F.
\end{equation*}
Moreover $\t(-\a,w)$, which is the unique affine function from $\Aa$
to $\RR$ whose linear part is ${-\a}$ and whose vanishing hyperplane
is $\Hh_{-a,w}$, vanishes on $x$.
Therefore it is equal to the map:
\begin{equation*}
y\mapsto-a(y)+a(x),
\end{equation*}
which implies that $\h_{-\a}(w)=a(x)$. 
According to {\bf B3} (see \S\ref{TroisGraces}), it follows that $w$
fixes $x$.

The subgroup $\U_{-a}(\F)$ is isomorphic to the additive group of
$\F$.
Thus, for any $r\in\RR$, the subgroup $\U_{-\a}(\F)_{r}$ corresponds
through this isomorphism to a nontrivial $\mathfrak{o}$-submodule of
$\F$, where $\mathfrak{o}$ denotes the ring of integers of $\F$ (see
\cite[Proposition 7.7]{Lan}).
Therefore there is a unique element $v\in\U_{-a}(\F)$ such that 
$w=v^2$ and $\h_{-\a}(v)=\h_{-\a}(w)$.
Thus $v\in\U_x$. 

The map $\U_{\a}(\F)\times\U_{\a}(\F)\to\G_\F$ defined by
$(u,u')\mapsto uwu'$ is injective and the intersection given by
(\ref{InterMU}) consists of a single element, namely $n$.
If we choose $u,u'\in\U_{\a}(\F)$ such that $uwu'=n$, then the
element:
\begin{equation}
\s(u')^{-1}w\s(u)^{-1}=\s(n)^{-1}
\end{equation}
is contained in the intersection (\ref{InterMU}).
Hence $\s(n)^{-1}$ is equal to $n$, and the unicity property implies
that $u'=\s(u)^{-1}$.
Moreover, according to \cite[Lemma 7.4(ii)]{Lan}, the real numbers 
$\h_{\a}(u)$ and $\h_{\a}(\s(u))$ are both equal to $-\h_{-\a}(w)$.
This implies that $u$ and $\s(u)$ are contained in $\U_x$.
Since $v$ is $\s$-anti-invariant and $w=v^{2}$, we get the expected
result with $c=(uv)^{-1}$.
\end{proof}

\begin{rema}
Note that $\s(c)\in\U_x$. 
Indeed we have $\s(v)=v^{-1}\in\U_{x}$ and $\s(u)\in\U_x$.
Hence $n=c^{-1}\s(c)\in\N_\F\cap\U_{\O}$, which is contained in
$\N_\O$ with $\O=\{x,\s(x)\}$. 
\end{rema}

\subsection{}
\label{Parsec}

Let $\Dd\G$ denote the derived subgroup of $\G$, and recall that $\C$
denotes the connected centre of $\G$.
This latter subgroup is a $\F$-split torus of $\G$.

\begin{lemm}
\label{Marguerites}
Let $\T$ be a $\F$-split torus of $\G$.
\begin{itemize}
\item[(i)]
There is a $\F$-subtorus $\T'$ of $\C$ such that the groups 
$\T\cdot\Dd\G$ and $\T'\cdot\Dd\G$ are equal.
\item[(ii)]
If $\T$ is $(\s,\F)$-split, then any $\T'$ satisfying (i) is
$(\s,\F)$-split.
\item[(iii)]
Assume that $\Dd\G$ is contained in $\H$ and $\T$ is $(\s,\F)$-split.
Then any $\T'$ satisfying (i) is $(\s,\F)$-split and has the same
dimension as $\T$. 
\end{itemize}
\end{lemm}

\begin{proof}
We set $\tilde\G=\G/\Dd\G$ and, for any $\F$-subgroup $\K$ of $\G$, we
denote by $\tilde\K$ the image of $\K$ in $\tilde\G$.
According to \cite[Proposition 14.2]{Borel}, the group $\G$ is the almost 
direct product of $\C$ and $\Dd\G$, which means that $\G$ is equal to
the product $\C\cdot\Dd\G$ and that the intersection $\C\cap\Dd\G$ is
finite. 
This implies that $\tilde\C=\tilde\G$.
Let $f$ denote the $\F$-rational map $\C\to\tilde\C$.
It is surjective with finite kernel.
Hence $\tilde\G$ is a $\F$-split torus, and we denote by $\tilde\s$
the involutive $\F$-automorphisme of $\tilde\G$ induced by $\s$.
We now prove the lemma in three steps.
\begin{itemize}
\item[(i)]
According to \cite[Proposition 8.2(c)]{Borel}, the neutral component
of the inverse image $f^{-1}(\tilde\T)$ is a $\F$-split subtorus of
$\C$ which we denote by $\T'$.
It has finite index in $f^{-1}(\tilde\T)$.
The image $f(\T')$ is then a subtorus of finite index in the connected
group $\tilde\T$, so that $\tilde\T'=\tilde\T$.
\item[(ii)]
Now assume that $\T$ is $(\s,\F)$-split, and let $\T'$ satisfy (i).
Let us consider the map $t\mapsto t\s(t)$ from $\T'$ to itself. 
As $\tilde\T'=\tilde\T$ is a $(\tilde\s,\F)$-split torus, the image of
this map is a connected $\F$-subgroup contained in the kernel of $f$,
which is finite.
\item[(iii)]
Assume that $\Dd\G$ is contained in $\H$ and $\T$ is $(\s,\F)$-split.
Then the map $\T\to\tilde{\T}$ has finite kernel, which implies that
$\T$ and $\tilde{\T}$ have the same dimension.
Now let $\T'$ satisfy (i).
According to (ii), such a torus is $(\s,\F)$-split, and it has the
same dimension as $\tilde\T'=\tilde\T$.
\end{itemize}
This ends the proof of Lemma \ref{Marguerites}.
\end{proof}

\subsection{}

Let $\Bb$ denote the building of $\G$ over $\F$.

\begin{prop}
\label{Rubempre}
Let $x$ be a special point of $\Bb$.
There exists a $\s$-stable max\-imal $\F$-split torus $\SS$ of $\G$
such that the apartment corresponding to $\SS$ cont\-ains $x$ and such
that $\SS^-$ is a maximal $(\s,\F)$-split torus of $\G$.
\end{prop}

\begin{rema}
In \S\ref{CounterStrike} we give an example of a {\it non split}
$\F$-group $\G$ such that Proposition \ref{Rubempre} does not hold.
\end{rema}

\begin{proof}
Let $\Aa$ be a $\s$-stable apartment containing $x$ (see Proposition
\ref{Raince}) and let $\SS$ be the corresponding maximal $\F$-split
torus of $\G$.
Assume that $\Aa$ has been chosen such that the dimension of the
$(\s,\F)$-split torus $\SS^-$ is maximal.
If it is a maximal $(\s,\F)$-split torus of $\G$, then Proposition
\ref{Rubempre} is proved. 
Assume that this is not the case, and let $\A$ be a maximal
$(\s,\F)$-split torus of $\G$ containing $\SS^{-}$.
The dim\-ension of $\A$ is greater than $\dim \SS^{-}$ (if not, the 
containment $\SS^{-}\subset\A$ would imply that $\SS^{-}=\A$).
If we get a contradiction, the proposition will be proved.

Let $\G'$ be the neutral component of the centralizer of $\SS^{-}$ in
$\G$.
It is a $\F$-split connected reductive subgroup of $\G$ containing
$\SS$ and $\A$, which is naturally endowed with a nontrivial action
of $\s$.
Let $\C'$ denote the connected center of $\G'$.

\begin{lemm}
\label{Petunias}
There is $a\in\Phi(\G',\SS)$ such that the corresponding root subgroup
$\U'_a$ is not contained in $\H$, and such a root is $\s$-invariant.
\end{lemm}

\begin{proof}
Assume that $\U'_a\subset\H$ for each root $a\in\Phi(\G',\SS)$.
Thus the derived subgroup $\Dd\G'$, which is generated by the
$\U'_a$ for $a\in\Phi(\G',\SS)$, is contained in $\H$
(see \cite[Theorem 27.5(e)]{Humphreys}).
According to Lemma \ref{Marguerites}(iii), there exists a
$(\s,\F)$-subtorus $\A'$ of $\C'$ such that
$\A\cdot\Dd\G'=\A'\cdot\Dd\G'$ and $\dim(\A)=\dim(\A')$.

The subgroup generated by $\C'$ and $\SS$ is a $\F$-torus of $\G'$.
As $\G'$ is $\F$-split, $\SS$ is a maximal torus of $\G'$, hence it 
contains $\C'$.
Therefore $\SS^{-}$ contains $\A'$ which has the same dimension as
$\A$, and this dimension is greater than $\dim\SS^{-}$.
This gives us a contra\-diction.

Now let $\a$ be a root in $\Phi(\G',\SS)$ such that $\U'_a$ is not
contained in $\H$.
The root $\a$ and its conjugate $\a\circ\s$ coincide on $\SS^+$ and
are both trivial on $\SS^-$. 
As $\SS$ is the almost direct product of $\SS^+$ and $\SS^-$ (see 
\S\ref{AlmostTrue}), they are equal. 
Therefore $\a$ is $\s$-invariant.
This ends the proof of Lemma \ref{Petunias}.
\end{proof}

Let $a\in\Phi(\G',\SS)$ as in Lemma \ref{Petunias}.
If we think to $\a$ as a root in $\Phi(\G,\SS)$, the root subgroup
$\U_\a$ is $\s$-stable and is not contained in $\H$.
Moreover, we have the following result.

\begin{lemm}
$\U_\a(\F)$ is contained in $\{g\in\G_\F\ |\ \s(g)=g^{-1}\}$.
\end{lemm}

\begin{proof}
As $\G$ is $\F$-split, $\U_\a$ is $\F$-isomorphic to the additive
group.
Thus the action of $\s$ on $\U_\a(\F)$ corresponds to an involutive
automorphism of the $\F$-algebra $\F[t]$.
It has the form $t\mapsto\l t$ for some $\l\in\mult\F$ with
$\l^{2}=1$. 
As $\U_\a$ is not contained in $\H$, we have $\l=-1$.
This gives the expected result.
\end{proof}

According to Lemma \ref{Almani}, there are $n\in\N_\F$ and $c\in\U_x$
such that $n=c^{-1}\s(c)$ and $\nu(n)$ is the affine reflection of
$\Aa$ which let $x$ invariant and whose linear part is $s_a$.
For any $t\in\SS$, note that we have:
\begin{eqnarray*}
\s(ctc^{-1})&=&cn\s(t)n^{-1}c^{-1}\\
&=&cs_\a(\s(t))c^{-1}.
\end{eqnarray*}
Let $\Aa'$ denote the apartment $c\cdot\Aa$ and let $\SS'={}^{c}\SS$
be the corresponding maximal $\F$-split torus of $\G$.
Then $\Aa'$ contains $x$ and is $\s$-stable. 
Moreover, as the root $a$ is trivial on $\SS^{-}$ and $s_\a$ fixes the
kernel of $\a$ pointwise, the conjugate ${}^c(\SS^{-})$ is a
$(\s,\F)$-split subtorus of $\SS'$.
Thus $\SS'^-$ has dimension not smaller than $\dim\SS^-$.

Now let $\SS_\a$ denote the maximal $\F$-split torus in the set of all
$t\in\SS$ such that $s_\a(t)=t^{-1}$.
As $\a$ is $\s$-invariant, such a torus is $\s$-stable.
Moreover, it is one dimensional and its intersection with $\Ker(\a)$
is finite.
Therefore the conjugate ${}^c\SS_\a$ is a nontrivial $(\s,\F)$-split
subtorus of $\SS'$ which is not contained in ${}^c(\SS^{-})$.
Thus the dimension of $\SS'^-$, which contains
${}^{c}(\SS_{\a}\SS^{-})$, is greater than $\dim\SS^{-}$, which
contradicts the maximality property of $\Aa$.
This ends the proof of Proposition \ref{Rubempre}.
\end{proof}

\subsection{}

Let $\A$ be a max\-imal $(\s,\F)$-split torus of $\G$ and $\SS$ a
$\s$-stable maximal $\F$-split torus of $\G$ containing $\A$.
Let $\{\A^{j}\ |\ j\in\J\}$ be a set of representatives of the
$\H_\F$-conjugacy classes of maximal $(\s,\F)$-split tori in $\G$.
Let $x$ be a special point of the building and let $\K$ be its
stabilizer in $\G_\F$.

\begin{theo}
\label{SplendeursEtMiseres}
For $j\in\J$, let $y_j\in\G_{\F}$ such that ${}^{y_{j}}\A=\A^{j}$.
We have:
\begin{equation*}
\G_{\F}=\bigcup\limits_{j\in\J}\H_\F y_{j}\SS_\F\K.
\end{equation*}
\end{theo}

\begin{proof}
We fix $g\in\G_\F$.
According to Proposition \ref{Rubempre}, there is a $\s$-stable
max\-imal $\F$-split torus $\SS'$ of $\G$ such that the apartment
corresponding to it cont\-ains $g\cdot x$ and such that $\SS'^-$ 
is a maximal $(\s,\F)$-split torus of $\G$.
Let $j\in\J$ be such that $\SS'^-$ is $\H_\F$-conjugate to $\A^j$.
According to \cite[Lemma 2.2]{HH}, there is $h\in\H_\F$ such that
$\SS'={}^{hy_j}\SS$.
Hence $g\cdot x$ is contained in $hy_j\cdot\Aa$.
According to Property (2) of \S\ref{TroisParques}, there exists
$n\in\N_{\F}$ such that $g\cdot x=hy_jn\cdot x$. 

Therefore $\G_\F$ is the union of the $\H_\F y_j\N_\F\K$ for
$j\in\J$. 
As $x$ is special, we have $\N_\F\K=\SS_\F\K$ and we get the expected
result. 
\end{proof}


\section{Examples}

Let $\F$ be a non-Archimedean locally compact field of residue
characteristic different from $2$.
Let $\o$ be its ring of integers and $\p$ its maximal ideal. 

\subsection{}

Here we consider the connected reductive $\F$-group $\G=\GL_{n}$,
endowed with the $\F$-involution $\s:g\mapsto{}^tg^{-1}$, where 
${}^{t}g$ denotes the transpose of $g\in\G$.
We set $\K=\GL_n(\o)$ and $\H=\G^\s$, which is an orthogonal group,
and we denote by $\SS$ the diagonal torus of $\G$.

We start with the following lemma.

\begin{lemm}
\label{lem2}
Let $\V$ be a finite dimensional $\F$-vector space and $\B$ a
symmetric bilinear form on $\V$.
Then any free $\o$-submodule of finite rank of $\V$ has a basis which
is orthogonal relative to $\B$.
\end{lemm}

\begin{proof}
Let $\La$ be a free $\o$-submodule of finite rank of $\V$. 
The proof goes by induction on the rank of $\La$. 
If $\B$ is null, then the result is trivial.
If not, we denote by $\B_{\La}$ the restriction of $\B$ to
$\La\times\La$.
Its image is of the form $\p^m$ for some integer $m\in\ZZ$.
If $\w$ is a uniformizer of $\F$, then the form
$\B_{\La}^0=\w^{-m}\B_{\La}$ has image $\o$ on $\La\times\La$.
Therefore, it defines a non trivial bilinear form:
\begin{equation*}
\bar\B_{\La}^0:\La/\p\La\times\La/\p\La\to\o/\p.
\end{equation*}
Let $e\in\La$ be a vector whose reduction mod. $\p$ is not isotropic
relative to $\bar\B_{\La}^0$, which means that $\B_{\La}^0(e,e)$ is a
unit of $\o$.
Then $\La$ is the direct sum of $\o e$ and $\La\cap\F e^{\bot}$,
where $\F e^{\bot}$ denotes the orthogonal of $\F e$ in $\V$.
Indeed, it follows from the decomposition:
\begin{equation*}
x=\frac{\B(e,x)}{\B(e,e)}e+\Big(x-\frac{\B(e,x)}{\B(e,e)}e\Big)
\end{equation*}
for any $x\in\La$.
As $\La\cap\F e^{\bot}$ is a free $\o$-submodule of finite rank of
$\V$ whose rank is smaller than the rank of $\La$, we conclude by
induction.
\end{proof}

We introduce the set $\Ee$ of all $g\in\G_\F$ such that
${}^tgg\in\SS_\F$ (compare \S\ref{IciOEstDefini}).
We have the following decomposition of $\G_\F$, which is more 
precise than the one given by Proposition \ref{Euripide}.

\begin{prop}
\label{decGSK}
We have $\G_\F=\Ee\K$.
\end{prop}

\begin{proof}
We make $\G_\F$ act on the quotient $\G_\F/\K$, which can be
identified to the set of all $\o$-lattices (that is, cocompact free
$\o$-submodules) of the $\F$-vector space $\V=\F^n$. 
Let $\B$ denote the symmetric bilinear form on $\V$ making the
canonical basis of $\V$ into an orthonormal basis.
According to Lemma \ref{lem2}, for any $g\in\G_\F$, the
$\o$-lattice $\La$ corresponding to the class $g\K$ has a basis which
is orthogonal relative to $\B$.
This means that there exists $u\in\K$ such that the element 
$g'=gu^{-1}\in g\K$ maps the canonical basis of $\V$ to an orthogonal
basis of $\La$. 
Therefore we have $g'\in\Ee$, thus $g\in\Ee\K$.
\end{proof}

We now investigate the maximal $(\s,\F)$-split tori of $\G$.
Note that $\SS$ is a maximal $(\s,\F)$-split torus of $\G$.

\begin{prop}
\label{Verneuil}
The map $g\mapsto{}^g\SS$ induces a bijection between
$(\H_\F,\N_\F)$-double cosets of $\Ee$ and $\H_\F$-conjugacy classes
of maximal $(\s,\F)$-split tori of $\G$.
\end{prop}

\begin{proof}
One immediately checks that this map is well defined and injective. 
For $g\in\G_\F$, the conjugate ${}^g\SS$ is a maximal $(\s,\F)$-split
torus of $\G$ if anf only if $g^{-1}\s(g)\in\SS_\F$, which amounts to
saying that $g\in\Ee$ and proves surjectivity.
\end{proof}

Let $\Qq$ denote the set of all equivalence classes of non
degenerate quadratic forms on $\F^n$.
For $a={\rm diag}(a_1,\ldots,a_n)\in\SS_\F$ we denote by $\Q_a$ the
diagonal quad\-ratic form $a_1\X_{1}^{2}+\ldots+a_n\X_{n}^{2}$.
Note that the map $a\mapsto\Q_a$ induces a surjective map from $\SS_\F$
to $\Qq$.

\begin{prop}
\label{Sprint1}
\begin{itemize}
\item[(i)]
The map $g\mapsto{}^tgg$ induces an injection $\iota$ from the set of 
$(\H_\F,\N_\F)$-double cosets of $\Ee$ to $\H^1(\N_\F)$. 
\item[(ii)]
For $a\in\SS_\F$, the class of $a$ in $\H^1(\N_\F)$ is in the image
of $\iota$ if and only if $\Q_a\sim\X_{1}^{2}+\ldots+\X_{n}^{2}$.
\end{itemize}
\end{prop}

\begin{proof}
We have an exact sequence:
\begin{equation*}
\H_\F\to\H^0(\G_\F/\N_\F)\to\H^1(\N_\F)\to\H^1(\G_\F),
\end{equation*}
where the map from $\H^0(\G_\F/\N_\F)$ to $\H^1(\N_\F)$
is induced by $g\mapsto{}^tgg$. 
As the set of $(\H_\F,\N_\F)$-double cosets of $\Ee$ is a subset of
$\H_\F\backslash\H^0(\G_\F/\N_\F)$, we get (i).
To get (ii), it is enough to remark that $\H^1(\G_\F)$ canonically
identifies with $\Qq$.
\end{proof}

\begin{rema}
Recall (see \cite[IV \S2.3]{Serre}) that for $a,b\in\SS_\F$, the
nondegenerate quadratic forms $\Q_a,\Q_b$ are equivalent if, and only
if they have the same discriminant and the same Hasse invariant.
\end{rema}

\begin{prop}
\label{Lollo}
Let $\{a^j\ |\ j\in\J\}\subset\SS_\F$ form a set of representatives of
$\Im(\iota)$ in $\H^1(\N_\F)$. 
For $j\in\J$, we choose $y_j\in\Ee$ such that ${}^ty_jy_j=a^j$.
Then:
\begin{equation}
\label{Binou}
\G_\F=\bigcup\limits_{j\in\J}\H_\F y_j\SS_\F\K.
\end{equation}
\end{prop}

\begin{proof}
Propositions \ref{decGSK} and \ref{Verneuil} imply that $\G_\F$ is the
union of the components $\H_\F y_j\N_\F\K$ for $j\in\J$.
As $\N_\F\K=\SS_\F\K$ we get the expected result.
\end{proof}

\begin{exem}
In the case where $n=2$, we give an explicit description of
$\Im(\iota)$.
Let $\w$ denote a uniformizer of $\o$ and $\xi\in\mult\o$ a non square
unit of $\o$, so that $\{1,\xi,\w,\xi\w\}$ is a set of
rep\-resentatives of $\mult\F$ modulo $\F^{\times 2}$.
The set of elements of $\mult\F$ which are represented by the
quadratic form $\Q_1=\X^2+\Y^2$ depends on the image of $p$ in
$\ZZ/4\ZZ$.
If $p\equiv 1\ {\rm mod.}\ 4$, all elements of $\mult\F$ are
represented by $\Q_1$.
If $p\equiv 3\ {\rm mod.}\ 4$, an element of $\mult\F$ is
represented by $\Q_1$ if and only if its normalized valuation if
even. 
We set:
\begin{equation*}
\J=
\left\{
\begin{array}{ll}
\{1,\xi,\w,\xi\w\}\quad & \text{if}\ p\equiv 1\ {\rm mod.}\ 4, \\
\{1,\xi\} & \text{if}\ p\equiv 3\ {\rm mod.}\ 4.
\end{array}
\right.
\end{equation*}
For each $j\in\J$, set $a^{j}={\rm diag}(j,j)$.
Then the elements $a^{j}$ form a set of rep\-resentatives of
$\Im(\iota)$ in $\H^1(\N_\F)$.
\end{exem}

\subsection{}
\label{Restrike}

\def\Perm{\mathfrak{S}}
\def\E{\F'}

In this paragraph we consider the connected reductive $\F$-group
$\G=\Res_{\E/\F}\GL_{n}$, where $\E$ is a quadratic
extension of $\F$, endowed with the involutive $\F$-auto\-morphism
$\s$ of $\G$ induced by the nontrivial element of $\Gal(\E/\F)$.

We set $\H=\G^{\s}$, so that we have $\G_\F=\GL_n(\E)$ and
$\H_\F=\GL_n(\F)$. 
We denote by $\SS$ the diagonal torus of $\G$ and by $\K$ the maximal
compact subgroup $\GL_n(\o')$ of $\G_\F$, where $\o'$ denotes the
ring of integers of $\E$.
Note that $\SS$ is $\s$-invariant, that is $\SS=\SS^+$.

As usual, $\N$ (resp. $\Z$) denotes the normalizer (resp. the
centralizer) of $\SS$ in $\G$.
Let $\Perm_n$ denote the group of permutation matrices in $\G_\F$, so
that $\N_\F$ is the semidirect product of $\Perm_n$ by $\Z_\F$.
Note that $\SS_\F$ (resp. $\Z_\F$) is the subgroup of all diagonal
matrices of $\G_\F$ with entries in $\F$ (resp. in $\E$).

\begin{lemm}
\label{Ploum}
$\H^1(\N_\F)$ can be identified with the set of conjugacy classes of
elements of $\Perm_n$ of order $1$ or $2$.
\end{lemm}

\begin{proof}
According to Hilbert's Theorem 90, the group $\H^1(\Z_\F)$ is trivial.
Therefore we have an exact sequence:
\begin{equation}
\label{hyppolite}
1\f\H^1(\N_\F)\f\H^1(\N_\F/\Z_\F).
\end{equation}
As $\s$ acts trivially on $\N_\F/\Z_\F\simeq\Perm_n$, the set
$\H^1(\N_\F/\Z_\F)$ can be identified to the set of
$\Perm_n$-conjugacy classes of $\Hom(\ZZ/2\ZZ,\Perm_n)$, that is, 
to the set of conjugacy classes of elements of $\Perm_n$ of order
$\<2$.
This proves that $\H^1(\N_\F)$ can be naturally embedded in the set of
conjugacy classes of elements of $\Perm_n$ of order $\<2$.

Now two elements $w,w'\in\Perm_n$ define the same class in
$\H^1(\N_\F)$ if and only if they are conjugate in $\Perm_n$, 
thus if and only if $w\Z_\F$ and $w'\Z_\F$ define the same class in 
$\H^1(\N_\F/\Z_\F)$. 
Therefore (\ref{hyppolite}) is a bijection.
\end{proof}

\begin{prop}
\label{Eulalie}
\begin{itemize}
\item[(i)]
The number of $\H_\F$-conjugacy classes of $\s$-stable maximal
$\F$-split tori in $\G_\F$ is $[n/2]+1$.
\item[(ii)]
There is a unique $\H_\F$-conjugacy class of maximal $(\s,\F)$-split
tori in $\G_\F$.
\end{itemize}
\end{prop}

\begin{proof}
\begin{itemize}
\item[(i)]
Let $\Oo$ denote the set of all $g\in\G_\F$ such that
$g^{-1}\s(g)\in\N_\F$. 
Then the map $g\mapsto{}^g\SS$ defines an injective map from the set
of $(\H_\F,\N_\F)$-double cosets of $\Oo$ to $\H^1(\N_\F)$.
Therefore we are reduced to proving that this map is surjective, and
(i) will follow from Lemma \ref{Ploum}

For $n=2$, let $\tau$ denote the nontrivial element of $\Perm_2$ and 
choose an element $a\in\E$ which is not in $\F$.
Then the element:
\begin{equation}
\label{DEFu}
u=
\begin{pmatrix}
a &\s(a)\\
1 &1\\
\end{pmatrix}
\in\GL_2(\E)
\end{equation}
satisfies the relation $u^{-1}\s(u)=\tau$.
For an arbitrary integer $n\>2$, let $w\in\Perm_n$ have order $\<2$.
Then there is an integer $0\<i\<[n/2]$ such that $w$ is 
conjugate to the element: 
\begin{equation*}
\tau_i={\rm diag}(\tau,\ldots,\tau,1,\ldots,1)\in\GL_n(\E),
\end{equation*}
where $\tau\in\GL_2(\E)$ appears $i$ times and $1\in\GL_1(\E)$ appears
$n-2i$ times. 
Thus the matrice:
\begin{equation}
\label{DEFui}
u_i={\rm diag}(u,\ldots,u,1,\ldots,1)\in\GL_n(\E)
\end{equation}
satisfies the relation $u_i^{-1}\s(u_i)=\tau_i$.
Therefore any cocycle in $\N_\F$ is $\G_\F$-cohomo\-logous to the
neutral element $1\in\G_\F$, which proves (i).
\item[(ii)]
For any $0\<i\<[n/2]$, the dimension of the $(\s,\F)$-split
torus $({}^{u_i}\SS)^-$ is equal to $i$. 
According to (i), the map:
\begin{equation*}
\H_\F g\N_\F\mapsto\text{class of } g^{-1}\s(g) \text{ in } \H^1(\N_\F)
\end{equation*}
is a bijection from the set of $(\H_\F,\N_\F)$-double cosets of $\Oo$
to $\H^1(\N_\F)$, and the elements of this latter set are the classes
of the $\tau_i$ for $0\<i\<[n/2]$.
This gives us the expected result.
\end{itemize}
This ends the proof of Proposition \ref{Eulalie}.
\end{proof}

\begin{prop}
For $0\<i\<[n/2]$, let $u_i$ denote the element defined by
(\ref{DEFu}) and (\ref{DEFui}).
Then we have:
\begin{equation*}
\G_\F=\bigcup_{i=0}^{[n/2]}
\H_\F u_i\Z_\F\K.
\end{equation*}
\end{prop}

\begin{proof}
According to the proof of Proposition \ref{Eulalie}, the set $\Oo$ is
the union of the double cosets $\H_\F u_i\N_\F$ with $0\<i\<[n/2]$.
The result then follows from Proposition \ref{Euripide} and from the
fact that $\N_\F\K=\Z_\F\K$.
\end{proof}

\subsection{}
\label{CounterStrike}
\def\AA{\mathfrak{A}}
\def\FF{{\rm F}}

Here we give an example (due to Bertrand Lemaire) of a non-split
$\F$-group such that Proposition \ref{Rubempre} does not hold. 
We set $\G=\Res_{\E/\F}\GL_{2}$, where $\E$ is now a 
{\it totally ramified} quadratic extension of $\F$.
The $\F$-involution $\s$ is still induced by the nontrivial element of
$\Gal(\E/\F)$ and we set $\H=\GL_{2}$.
Let $\Bb'$ (resp. $\Bb$) denote the building of $\G$ (resp. $\H$) over
$\F$.

In \cite{BT3}, Bruhat and Tits give a description of the
faces of $\Bb$ in terms of hereditary $\o$-orders of $\M_{2}(\F)$.
More precisely, there is a bijective correspondence:
\begin{equation}
\label{IdeeNoire}
\FF\mapsto\Mm_\FF
\end{equation}
between the faces of $\Bb$ and the hereditary $\o$-orders of
$\M_{2}(\F)$, such that the stabilizer of $\FF$ in $\GL_2(\F)$ in the
normalizer of $\Mm_\FF$ in $\GL_2(\F)$.
For $x\in\Bb$, we will denote by $\Mm_x$ the hereditary order
corresponding to the face of $\Bb$ which contains $x$.
Of course, we have a similar correspondence between faces of
$\Bb'$ and hereditary $\o'$-orders of $\M_{2}(\F')$.
Moreover, as $\E$ is tamely ramified over $\F$, there is a bijective
correspondence $j$ from the set $\Bb'^\s$ of $\s$-invariant points
of $\Bb'$ to $\Bb$ such that, for any $x\in\Bb'^\s$, we have:
\begin{equation*}
\Mm_{j(x)}=\Mm_x\cap\M_{2}(\F).
\end{equation*}

Let $q$ denote the cardinal of the residue field of $\F$.
As $\E$ is totally ramified over $\F$, any vertex of $\Bb$
(resp. $\Bb'$) has exactly $q+1$ neighbours in $\Bb$ (resp. in
$\Bb'$). 
Let $x$ be a $\s$-invariant point of $\Bb'$.
Recall that, according to Proposition \ref{Raince}, it is contained in
a $\s$-stable apartment. 
\begin{enumerate}
\item
If $j(x)$ is in a chamber of $\Bb$, then $x$
has $q+1$ neighbours in $\Bb'$ but only two $\s$-invariant ones.
Thus $x$ has non-$\s$-invariant neighbours.
\item
If $j(x)$ is a vertex of $\Bb$, then $x$
has $q+1$ neighbours in $\Bb'$ as in $\Bb$.
Thus any neighbour of $x$ in $\Bb'$ is $\s$-invariant, which implies
that any $\s$-stable apartment containing $x$ is $\s$-invariant.
For instance, this is the case of the vertex $x$ corresponding to the
$\o'$-order $\M_{2}(\o')$, because its image $j(x)$ corresponds to the
maximal $\o$-order $\M_{2}(\o')\cap\M_{2}(\F)=\M_{2}(\o)$.
Such a special point does not satisfy Proposition \ref{Rubempre}.
\end{enumerate}


\providecommand{\bysame}{\leavevmode ---\ }
\providecommand{\og}{``}
\providecommand{\fg}{''}
\providecommand{\smfandname}{\&}
\providecommand{\smfedsname}{\'eds.}
\providecommand{\smfedname}{\'ed.}
\providecommand{\smfmastersthesisname}{M\'emoire}
\providecommand{\smfphdthesisname}{Th\`ese}

\end{document}